\documentclass[11pt,reqno]{amsproc}

\title[Inviscid limit of Navier-Stokes for bounded velocity fields]{Remarks on the inviscid limit for the Navier-Stokes equations\\ for uniformly bounded velocity fields}

\author[P. Constantin]{Peter Constantin}
\address{Department of Mathematics, Princeton University, Princeton, NJ 08544}
\email{const@math.princeton.edu}

\author[T. Elgindi]{Tarek Elgindi}
\address{Department of Mathematics, Princeton University, Princeton, NJ 08544}
\email{tme2@math.princeton.edu}

\author[M. Ignatova]{Mihaela Ignatova}
\address{Department of Mathematics, Princeton University, Princeton, NJ 08544}
\email{ignatova@math.princeton.edu}

\author[V. Vicol]{Vlad Vicol}
\address{Department of Mathematics, Princeton University, Princeton, NJ 08544}
\email{vvicol@math.princeton.edu}

\usepackage[margin=1.25in]{geometry}
\usepackage{amsmath, amsthm, amssymb}
\usepackage{times}
\usepackage{graphicx}

\usepackage{color}
\usepackage[colorlinks=true, pdfstartview=FitV, linkcolor=blue, citecolor=blue, urlcolor=blue]{hyperref}

\theoremstyle{plain}
\newtheorem{theorem}{Theorem}[section]

\theoremstyle{definition}
\newtheorem{remark}[theorem]{Remark}

\numberwithin{equation}{section}

\newcommand{\blu}[1]{#1}

\newcommand{\be}{\begin{equation}}
\newcommand{\eeq}{\end{equation}}

\def\RR{{\mathbb R}}
\def\HH{{\mathbb H}}

\def\OO{{\mathcal O}}

\def\eps{{\varepsilon}}

\def\phi{\varphi}

\def\erf{ \mathop{\rm erf} \nolimits}
\def\erfc{ \mathop{\rm erfc} \nolimits}

\def\uns{{u}^{\mbox{\tiny NS}}}
\def\CNS{{C}_{\mbox{\tiny NS}}}
\def\CE{{C}_{\mbox{\tiny E}}}
\def\pns{{p}^{\mbox{\tiny NS}}}
\def\vortns{{\omega}^{\mbox{\tiny NS}}}
\def\vorte{{\omega}^{\mbox{\tiny E}}}
\def\ue{{u}^{\mbox{\tiny E}}}
\def\UE{{U}^{\mbox{\tiny E}}}
\def\pe{{p}^{\mbox{\tiny E}}}
\def\uk{{u}^{\mbox{\tiny K}}}

\def\supp{\mathop{\rm supp} \nolimits}

\begin{document}

%%%%%%%%%%%%%%%%%%%%%%%% Abstract %%%%%%%%%%%%%%%%%%%%
\begin{abstract}
We consider the vanishing viscosity limit of the Navier-Stokes equations in a half space, with Dirichlet boundary conditions. We prove that the inviscid limit holds in the energy norm if the product of the components of the Navier-Stokes solutions \blu{are equicontinuous at $x_2=0$}. \blu{A sufficient condition for this to hold is that the tangential Navier-Stokes velocity remains uniformly bounded and has a uniformly integrable tangential gradient near the boundary.}
%\hfill \today
\end{abstract}

\maketitle

%%%%%%%%%%%%%%%%%%%%%%%%%% The Main Part %%%%%%%%%%%%%%%%%%%%%%%%%%%%%%%%%%

\section{Introduction}
Consider the 2D Navier-Stokes equations
\begin{align}
&\partial_t \uns + \uns \cdot \nabla \uns + \nabla \pns = \nu\Delta \uns \label{eq:NS:1} \\
&\nabla \cdot \uns = 0 \label{eq:NS:2} \\
& \uns_1|_{\partial \HH} = \uns_2|_{\partial \HH} = 0 \label{eq:NS:3}
\end{align}
with kinematic viscosity $\nu$, in the half space $\HH = \{ (x_1,x_2) \colon x_2 >0 \}$, and the Euler equations
\begin{align}
&\partial_t \ue + \ue \cdot \nabla \ue + \nabla \pe = 0 \label{eq:E:1}\\
&\nabla \cdot \ue = 0 \label{eq:E:2} \\
& \ue_2|_{\partial \HH} = 0 \label{eq:E:3}
\end{align}
with asymptotically matching initial conditions
\begin{align}
\lim_{\nu \to 0} \| \uns_0 - \ue_0\|_{L^2(\HH)} = 0.
\label{eq:IC:cond}
\end{align}
We  denote by
\begin{align*}	
\ue_1|_{\partial \HH}(x_1,t) = \UE(x_1,t)
\end{align*}
the trace on $\partial \HH$ of the Euler tangential flow. We omit $\nu$ in the notation for $\uns$. Throughout this paper we consider $0 < \nu \leq \nu_0$, and $0 \leq t\leq  T$, where $\nu_0$ is an arbitrary fixed kinematic viscosity, and $T$ is an arbitrary fixed time. 
We assume that the Euler initial datum is smooth,  $\ue_0 \in H^s(\HH)$ for some $s>2$, so that there exists an unique $H^s$ smooth solution $\ue$ of \eqref{eq:E:1}--\eqref{eq:E:3} on $[0,T]$.

This paper establishes sufficient conditions for the family of
Navier-Stokes solutions $\{\uns\}_{\nu \in (0,\nu_0]}$ to ensure that the inviscid limit holds in the energy norm:
\begin{align}
\lim_{\nu \to 0} \|\uns  - \ue \|_{L^\infty(0,T;L^2 (\HH))} = 0.
\label{eq:inviscid:limit}
\end{align}
Our main results are given in Theorems~\ref{thm:main}, \ref{cor:2}, and \ref{thm:1}.
\blu{The main assumptions are either on the equicontinuity of the product $\uns_1 \uns_2$ at $x_2=0$, or on the uniform boundedness of $\uns_1$ and the uniform integrability of $\partial_1 \uns_1$ near the boundary of the half space. These conditions do not assume a priori any particular scaling law of the boundary layer with respect to $\nu$.} The conditions imposed imply that the Lagrangian paths originating in a boundary layer, stay in a proportional boundary layer during the time interval considered. The physical interpretation of our result is that, as long as there is no separation of the boundary layer, the inviscid limit  is possible.

\subsection{Known finite time, inviscid limit results}
The question of whether  \eqref{eq:inviscid:limit} holds in the case of Dirichlet boundary conditions has a rich history. Kato proved in~\cite{Kato84b} that the inviscid limit holds in the energy norm if and only if 
\begin{align}
 \lim_{\nu \to 0} \nu \int_0^T \!\!\! \int_{|x_2|\leq C \nu} |\nabla \uns(x_1,x_2,t)|^2 dx_1 dx_2 dt = 0,
 \label{eq:Kato}
\end{align}
i.e. that the energy dissipation rate is vanishing in a thin, $\OO(\nu)$, layer near the boundary. Kato's criterion was revisited and sharpened by many authors. For instance, in~\cite{TemamWang97b} and \cite{Wang01} it is shown that the condition on the full gradient matrix $\nabla \uns$ may be replaced by a condition on the tangential gradient of the Navier-Stokes solution alone, at the cost of considering a thicker boundary layer, of size $\delta(\nu)$, where $\lim_{\nu \to 0} \delta(\nu) / \nu =0$. In~\cite{Kelliher07} it is shown that $\nu \|\nabla \uns \|_{L^2(|x_2|\leq C \nu)}^2$ may be replaced by $\nu^{-1} \| \uns\|\blu{^2}_{L^2(|x_2|\leq C \nu)}$ which has the same scaling in the Kato layer. In~\cite{Kelliher08} it is shown that \eqref{eq:inviscid:limit} is equivalent to the weak convergence of vorticities
\begin{align}
 \vortns \to \vorte - \ue_1 \, \mu_{\partial \HH} \quad \mbox{in} \quad (H^1(\HH))^*
 \label{eq:JK}
\end{align}
where $\mu_{\partial \HH} $ is the Dirac measure on $\partial \HH$, and $(H^1)^*$ is the dual space to $H^1$ (not $H^1_0$). In fact, it is shown in~\cite{BardosTiti07} that the weak convergence of vorticity on the boundary 
\begin{align}
 \nu \vortns \to 0 \quad \mbox{in} \quad {\mathcal D}'([0,T]\times \partial\HH)
 \label{eq:BT}
\end{align}
is equivalent to \eqref{eq:inviscid:limit} (see also~\cite{Kelliher08,ConstantinKukavicaVicol15} in the case of stronger convergence in \eqref{eq:BT}). 

The idea to introduce a boundary layer corrector like Kato's, which is not based on power series expansions, and to treat the remainders with energy estimates has proven to be very fruitful. See for instance: \cite{Masmoudi98} in the case of anisotropic viscosity; \cite{BardosSzekelyhidiWiedemann14,BardosTitiWiedemann12} in the context of weak-strong uniqueness; \cite{GuoNguyen14} for a steady flow on a moving plate; \cite{BardosNguyen14} for the compressible Navier-Stokes equations; \cite{Lopes2TitiZang14} for the vanishing $\alpha$ limit of the 2D Euler-$\alpha$ model.

There are three classes of functions for which there exist unconditional inviscid limit results, that is, theorems whereby conditions imposed solely on initial data guarantee that \eqref{eq:inviscid:limit} is true for a time interval independent of viscosity (but possibly depending on initial data). The first class is that of real analytic initial data in all space variables~\cite{SammartinoCaflisch98b}, the second is that of initial data with vorticity
supported at an $\OO(1)$ distance from the boundary~\cite{Maekawa14}, and the third class is data with certain symmetries or special restrictions~\cite{LopesMazzucatoLopes08,LopesMazzucatoLopesTaylor08,MazzucatoTaylor08,Kelliher09}).
It is worth noting that in these three cases the Prandtl expansion of the Navier-Stokes equation is valid in a boundary layer of thickness $\sqrt{\nu}$. Moreover, in all these results, the Kato criteria also hold~\cite{BardosTiti07,Kelliher14}.  However, to date, there is no robust connection between the well-posedness of the Prandtl equations, and the vanishing viscosity limit in the energy norm. 

It is known that for a class of initial conditions close to certain shear flows the Prandtl equations are ill-posed~\cite{GerardVaretDormy10,GuoNguyen11,GerardVaretNguyen12} and even that the Prandtl expansion is not valid~\cite{Grenier00,GrenierGuoNguyen14,GrenierGuoNguyen14b,GrenierGuoNguyen14c}. These results do not imply that the inviscid limit in the energy norm is invalid, but rather just that the Prandtl expansion does not describe the leading order behavior near the boundary.  It would be natural to expect that working in a function space for which the local existence of the Prandtl equations holds  (see, e.g.~\cite{Oleinik66,MasmoudiWong12a,AlexandreWangXuYang14}, \cite{SammartinoCaflisch98a,LombardoCannoneSammartino03,KukavicaVicol13a}, \cite{KukavicaMasmoudiVicolWong14}, \cite{GerardVaretMasmoudi13}), there is a greater chance for \eqref{eq:inviscid:limit} to be true. An instance of such a result is given in~\cite{ConstantinKukavicaVicol15}, where a one-sided Kato criterion in terms of the vorticity is obtained, connecting Oleinik's monotonicity assumption and the inviscid limit: if
\begin{align}
\lim_{\nu \to 0} \int_0^T \left\| \left( \UE(x_1,t) \left( \vortns(x_1,x_2,t) + \frac{\delta(\nu t)}{\nu t} \right) \right)_- \right\|_{L^2(|x_2| \leq \nu t /\delta(\nu t))}^2 dt = 0
\label{eq:CKV}
\end{align}
holds, where  $\int_0^T \delta(\nu t) dt \to 0$ as $\nu \to 0$, then \eqref{eq:inviscid:limit} holds. In particular, if there is no back-flow in the underlying Euler flow, $\UE \geq 0$, and the Navier-Stokes vorticity $\vortns$ is larger than $- \delta(\nu t)/ \nu$ (for instance if it is non-negative as in Oleinik's setting) in a boundary layer that is slightly thicker than Kato's, then the inviscid limit holds. 

In contrast to the works \eqref{eq:Kato}--\eqref{eq:CKV} mentioned above, the goal of this paper is to establish sufficient conditions for \eqref{eq:inviscid:limit} to hold, which do not rely on any assumptions concerning derivatives of the Navier-Stokes equations. Alternately, we establish conditions which require only $L^1$ uniform integrability of tangential derivatives near the boundary. Our proofs keep the idea of Kato of building an ad-hoc boundary layer corrector, but its scaling is dictated by the heat equation in $x_2$ (with Prandtl scaling). No explicit convergence rates are obtained with our assumptions. The main results of this paper are:

\subsection{Results}
\begin{theorem}\label{thm:main}
Assume  that the family
\begin{align}
\{ \uns_1 \uns_2\}_{\nu \in (0,\nu_0]} \quad \mbox{is equicontinuous at} \quad x_2 = 0. 
\label{eq:assume:2}
\end{align}
Then \eqref{eq:IC:cond} implies that the inviscid limit holds in the energy norm.
\end{theorem}

 Specifically, in view of the Dirichlet boundary condition \eqref{eq:NS:3}, by condition \eqref{eq:assume:2} we mean  that there exists a function
\begin{align}
0 \leq \gamma (x_1,t) \in L^1_{t,x_1}([0,T] \times \RR)
\label{eq:assume:2a}
\end{align}
with the property that for any $\eps>0$, there exists $\rho = \rho (\eps) >0$ such that
\begin{align}
|\uns_1(x_1,x_2,t) \uns_2(x_1,x_2,t)| \leq \eps \gamma(x_1,t),  \quad \mbox{for all} \quad  x_2 \in (0,\rho],
\label{eq:assume:2b}
\end{align}
and all $(t,x_1) \in [0,T] \times \RR$, uniformly in $\nu \in (0,\nu_0]$. 

\blu{
\begin{remark}
We note that condition \eqref{eq:assume:2} holds for the solution of the Stokes equation (the linearization of the Navier-Stokes equations around the trivial flow $(0,0)$) and for the solution of  the Oseen equations (the linearization of the Navier-Stokes equations around a stationary shear flow $(U(y),0)$), {\em if one considers sufficiently smooth and compatible initial datum}. Indeed, the equicontinuity of the family $\{u_1 u_2\}_{\nu \in (0,\nu_0]}$ at $\{ x_2 = 0\}$ follows (with an explicit rate) for instance from the bound
\begin{align*}
\left| u_1(x_1,x_2,t) u_2(x_1,x_2,t) \right| \leq |u_1(x_1,x_2,t)| \int_0^{x_2} |\partial_1 u_1(x_1,y,t)| dy 
\end{align*}
and the fact that both $u_1$ and $\partial_1^k u_1$ obey the same equation for any $k \geq 1$. That is, since $\partial_1$ is a tangential vector field to $\partial \HH$, so that it commutes with the Stokes operator, for initial datum that is compatible one may obtain the same good bounds for $\partial_1^k u_1$ for all $k \geq 0$ (one may use the Ukai formula directly for the Stokes equation, or the argument in~\cite{TemamWang95} for Oseen). 
\end{remark}
}

The \blu{integrability (uniform in $\nu$) of $\partial_1 \uns_1$ and the boundedness (uniform in $\nu$) of $\uns_1$} is thus related to condition~\eqref{eq:assume:2}:
\begin{theorem}
\label{cor:2}
\blu{Let $\CNS>0$ be a constant. Assume that for any $\eps>0$ there exists $\rho>0$, such that 
\begin{align}
\sup_{\nu \in (0,\nu_0]} \int_0^T \| \partial_1 \uns_1(t)\, {\bf 1}_{0<x_2 <\rho} \|_{L^1(\HH)}^2 dt \leq \eps,
 \label{eq:assume:G3} 
\end{align}
meaning that the tangential derivative of the tangential component of the Navier-Stokes flow is uniformly integrable near $\{x_2=0\}$, 
and that
\begin{align}
\sup_{\nu \in (0,\nu_0]} \int_0^T \| \uns_1(t) \, {\bf 1}_{0<x_2 <\rho} \|_{ L^\infty(\HH)}^2  dt \leq \CNS \nu_0,
\label{eq:assume:1}
\end{align}
meaning that the tangential component of the Navier-Stokes flow is bounded near $\{x_2=0\}$.
}
Then \eqref{eq:inviscid:limit} holds.
\end{theorem}

The quantity in condition \eqref{eq:assume:1} is natural to consider: it is scale invariant under the Navier-Stokes isotropic scaling, and it appears in three dimensions as well. 
A similar quantity was used in~\cite{BardosSzekelyhidiWiedemann14} to establish conditional weak-strong uniqueness of weak solutions in H\"older classes.

Condition \eqref{eq:assume:G3} requires that the family of measures
\[
 \mu_\nu (dx_1\,dx_2) = |\partial_1 \uns_1(t,x_1,x_2)| dx_1 \, dx_2
\]
is uniformly absolutely continuous at $x_2=0$ with values in $L^2(0,T)$, \blu{i.e. that they do not assign uniformly bounded from below mass to a boundary layer}. Note that $\partial_1 \uns_1$ vanishes identically on $\partial \HH$, which is not the case for the Navier-Stokes vorticity $\vortns= \partial_2 \uns_1 - \partial_1 \uns_2$, which is expected to develop a measure supported on the boundary of the domain in the inviscid limit~\cite{Kelliher08}. Thus, the vorticity is not expected to be uniformly integrable in $L^2_t L^1_x$. Therefore, in \eqref{eq:assume:G3} it is important that instead of a uniform integrability condition on $\vortns$ or equivalently $\partial_2 \uns_1$, we have only assumed a uniform integrability condition on $\partial_1 \uns_1$. Also, note that (uniform in $\nu$) higher integrability of the Navier-Stokes vorticity, such as $L^p$ for $p>2$ cannot hold unless $\UE \equiv 0$, as is shown in~\cite{Kelliher14}.  

\blu{
\begin{remark}
\label{rem:L2:x2}
From the proof of Theorem~\ref{cor:2} it follows that in \eqref{eq:assume:G3}--\eqref{eq:assume:1} we may replace (both) $L^1_{x_1,x_2}(\HH)$ by $L^2_{x_1} L^1_{x_2}(\HH)$ and $L^\infty_{x_1,x_2}(\HH)$ with $L^2_{x_1} L^\infty_{x_2}(\HH)$. That is, the boundedness of $\uns_1$ and the uniform integrability of $\partial_1 \uns_1$ is to be only checked with respect to the $x_2$ variable. With respect to $x_1$ we only need that these functions are square integrable. Again, either of these conditions may be checked directly for the Stokes and Oseen equations if the initial datum is sufficiently smooth (with respect to $x_1$) and compatible. 
\end{remark}
}

A similar result to the one in Theorem~\ref{cor:2}, has been obtained independently in~\cite{Anna15}, where the authors prove that if $\nabla \uns$ is uniformly in $\nu$ bounded in $L^\infty(0,T; L^1(\Omega))$, for a domain $\Omega$ such that the embedding $W^{1,1}(\Omega) \subset L^2$ is compact, then the vanishing viscosity limit holds in $L^\infty(0,T;L^2(\Omega))$. \blu{Note moreover, that in~\cite{Wang01}, the author establishes a Kato-type criterion for the inviscid limit to hold, under the assumptions that $\nu \int_{0}^T \int_{|x_2|\leq \delta(\nu)} |\partial_1 \uns_1|^2 dt \to 0$ as $\nu \to 0$ and $\nu / \delta(\nu) \to 0$ as $\nu\to 0$. The modification to condition \eqref{eq:assume:G3} of this paper, as described in Remark~\ref{rem:L2:x2}, is weaker than the condition of \cite{Wang01}, as may be seen by taking $\rho = \nu$ and applying the H\"older inequality to bound $L^1_{x_2}$ with $L^2_{x_2}$. However, we also need to impose condition \eqref{eq:assume:1}.}

We conclude the introduction by noting that a similar proof to that of Theorem~\ref{thm:main} yields:
\blu{
\begin{theorem}
\label{thm:1}
Let 
$\delta$ be an increasing non-negative function such that $\lim_{\nu \to 0} \delta(\nu t)= 0$ uniformly for $t\in [0,T]$, such that $(\nu t) \delta'(\nu t)/ \delta(\nu t)$ is uniformly bounded for $\nu \in (0,\nu_0]$ and $t \in [0,T]$, and such that
\begin{align}
\lim_{\nu \to 0}   \int_0^T \frac{\nu}{\delta(\nu t)} + \frac{\delta(\nu t)}{t^{1+c}} dt= 0
\label{eq:Wang:assumption}
\end{align}
holds for some $c>0$. We may for instance take $\delta(\nu t) = (\nu t)^a$ for any $a\in (0,1)$. 

Assume that 
\begin{align}
\lim_{\nu \to 0} \uns_1(x_1,\delta(\nu t)\, y, t) \uns_2(x_1,\delta(\nu t)\, y, t) = 0
\label{eq:assume:4}
\end{align}
holds for a.e. $(t,x_1,y) \in [0,T] \times \HH$, and that 
\begin{align}
\int_0^T \sup_{\nu \in (0,1]}   \| \uns_1 (t)  \uns_2(t) \|_{ L^\infty (\{x_2 \leq \delta(\nu t) (\log(1/\nu))^{1/2}\}) }   <\infty.
\label{eq:assume:3}
\end{align} 
Then \eqref{eq:inviscid:limit} holds. 
\end{theorem}
}

\subsection{Organization of the paper} In Section~\ref{sec:proof:1} we lay out the scheme of the proof for the above mentioned theorems, by identifying the principal error terms in the energy estimate for the corrected $\uns-\ue$ flow. In Section~\ref{sec:corrector} we build a caloric lift of the Euler boundary conditions, augmented by an $\OO(1)$ correction at unit scale. In Section~\ref{sec:proof:2} we conclude the proof of Theorem~\ref{thm:main}, in Section~\ref{sec:proof:3} we give the proof of Theorem~\ref{cor:2}, while in Section~\ref{sec:proof:4} we show why Theorem~\ref{thm:1} holds.

\section{Setup of the Proof of Theorem~\ref{thm:main}}
\label{sec:proof:1}

We consider a boundary layer corrector $\uk$ (to be constructed precisely later) which for now obeys three properties
\begin{align}
&\nabla \cdot \uk = 0 \label{eq:uk:1}\\
&\uk_1|_{\partial \HH} = - \UE \label{eq:uk:2}\\
&\uk_2|_{\partial \HH} = 0. \label{eq:uk:3}
\end{align}
The main difference between the corrector $\uk$ we consider, and the one considered in~\cite{Kato84b}, is its characteristic length scale: we let $\uk$ obey a Prandtl $\sqrt{\nu t}$ scaling \blu{(see also~\cite{TemamWang95,Gie14})}. Roughly speaking, $\uk_1$ is a lift of the Euler boundary condition which obeys the heat equation $(\partial_t - \nu \partial_{x_2 x_2}) \uk_1 = 0$ to leading order in $\nu$. In view of \eqref{eq:uk:1}--\eqref{eq:uk:3} we then obtain $\uk_2$ from $\uk_1$ as 
\begin{align}
\uk_2(x_1,x_2,t) = - \int_{0}^{x_2} \partial_1 \uk_1(x_1,y,t) dy.
\label{eq:uk:incompressible}
\end{align}

The function
\begin{align*}
v = \uns - \ue - \uk
\end{align*}
is divergence free
\begin{align*}
\nabla \cdot v = 0
\end{align*} 
and obeys Dirichlet boundary conditions
\begin{align*}
v|_{\partial \HH} = 0.
\end{align*}
The equation obeyed by $v$ is 
\begin{align}
&\partial_t v - \nu \Delta v + v \cdot \nabla \ue  + \uns \cdot \nabla v + \nabla q \notag\\
&\qquad = \nu \Delta \ue - \left( \partial_t \uk - \nu \Delta \uk + \uns \cdot \nabla \uk + \uk \cdot \nabla \ue \right)
\label{eq:v:1}
\end{align}
where $q = \pns - \pe$.  Multiplying \eqref{eq:v:1} with $v$ and integrating by parts, yields 
\begin{align}
\frac 12 \frac{d}{dt} \|v\|_{L^2}^2 + \nu \|\nabla v\|_{L^2}^2 
&\leq \|\nabla \ue\|_{L^\infty} \|v\|_{L^2}^2  + \nu \|\Delta \ue\|_{L^2} \|v\|_{L^2} \notag\\
&\quad + T_1 + T_2 +T_3 + T_4 + T_5 + T_6
\label{eq:v:2}
\end{align}
where we have denoted
\begin{align}
T_1 &= - \int_{\HH} (\partial_t \uk - \nu \Delta \uk) \cdot v \label{eq:T1:def}\\
T_2 &=  - \int_{\HH} (\uns \cdot \nabla \ue )\cdot \uk \\
T_3 &= - \int_{\HH}  (\uk \cdot \nabla \ue) \cdot v\\
T_4 &= - \int_{\HH}  \uns_1  \uns_2 \partial_1 \uk_2 \\
T_5 &= - \int_{\HH}  \left( (\uns_1)^2 - (\uns_2)^2 \right) \partial_1 \uk_1 \label{eq:T5:def} \\
T_6 &= - \int_{\HH} \uns_1  \uns_2 \partial_2 \uk_1 \label{eq:T6:def}
\end{align}
The corrector $\uk$ is designed to eliminate the contribution from $T_1$ to leading order in $\nu$.
In turn, this leads to $\|\uk\|_{L^2} + \| \partial_1 \uk\|_{L^2} \to 0$ as $\nu\to 0$, so that the terms $T_2, T_3$, $T_4$, \blu{and $T_5$} are harmless. Such  is the case if $\uk$ is localized in a layer near the boundary, which is vanishing as $\nu \to 0$.
\blu{Assumption \eqref{eq:assume:2} only comes into play in showing that $T_6$ is} bounded conveniently. The next section is devoted to the construction of an $\uk$ with these properties, and the conclusion of the proof is given in Section~\ref{sec:proof:2} below. 

Throughout the text we shall denote by $\CE$ any constant that depends on $\|\ue\|_{L^\infty(0,T;H^s(\HH))}$. Various other positive constants shall be denoted by $C$; these constants do not depend on $\nu$, but they are allowed to implicitly depend on the fixed length of the time interval $T$, and on the largest kinematic viscosity $\nu_0$.

\section{A pseudo-caloric lift of the boundary conditions}
\label{sec:corrector}
\blu{The type of corrector we construct here was also used for instance in~\cite{TemamWang95,Gie14} (see also references therein) to address the vanishing viscosity limit for the linear Stokes system with compatible, respectively non-compatible initial datum.}

\subsection{The tangential component of the lift $\uk$}
Let
\begin{align*}
z = z(x_2,t) = \frac{x_2}{\sqrt{ 4 \nu t}}
\end{align*}
be the self-similar variable for the heat equation in $x_2$, with viscosity $\nu$. 
Let  $\eta$ be a non-negative bump function such that 
\begin{align}
\supp(\eta ) \in [1,2] \quad \mbox{and} \quad \int_1^2 \eta(r) dr = \frac{1}{\sqrt{\pi}}
\label{eq:eta}
\end{align}
which in addition obeys that $|\eta '|_{L^\infty} + |\eta''|_{L^\infty} \leq C_\eta$, for some constant $C_\eta$.

We let $\uk_1$ consist of a caloric lift of the Euler boundary conditions, augmented with a localization factor at large values of $x_2$. We define
\begin{align}
\uk_1(x_1,x_2,t) = - \UE(x_1,t) \left(\erfc(z(x_2,t)) - \sqrt{4 \nu t} \, \eta(x_2) \right) \label{eq:uk:1:def}
\end{align}
where
\begin{align*}
\erfc(z) = 1- \erf(z) = \frac{2}{\sqrt{\pi}} \int_z^\infty \exp(-y^2) dy.
\end{align*}
The normalization of the mass of $\eta$ was chosen precisely so that 
\begin{align}
\int_0^\infty \uk_1(x_1,x_2,t) dx_2 
&= - \UE(x_1,t)  \int_0^\infty \left( \erfc(z(x_2,t)) - \sqrt{4 \nu t} \, \eta(x_2) \right) dx_2 \notag\\
&= - \UE(x_1,t)  \sqrt{4 \nu t} \left( \int_0^\infty \erfc(z) dz - \int_0^\infty \eta(x_2) dx_2 \right) \notag\\ 
&= 0.
\label{eq:uk:1:zero:mean}
\end{align}
Property \eqref{eq:uk:1:zero:mean} of $\uk_1$ allows the $\uk_2$ defined in \eqref{eq:uk:incompressible} (see also below) to decay sufficiently fast as $x_2 \to \infty$. This decay of $\uk_2$ will be used essentially later on in the proof.

Note that $\uk_1$ is pseudo-localized to scale $x_2 \approx \sqrt{4 \nu t}$. Indeed, we have that
\begin{align*}
\|  \erfc(z(x_2,t)) \|_{L^p_{x_2}(0,\infty)} 
&= (4\nu t)^{1/(2p)} \| 1 - \erf(z)\|_{L^p_z(0,\infty)} \notag\\
&\leq C (\nu t)^{1/(2p)}
\end{align*}
and 
\begin{align*}
\| \partial_{x_2} \erfc(z(x_2,t)) \|_{L^p_{x_2}(0,\infty)} 
&= (4\nu t)^{1/(2p) - 1/2} \|\partial_z  \erfc(z)\|_{L^p_z(0,\infty)} \notag\\
&\leq C (\nu t)^{1/(2p)-1/2}
\end{align*}
for all $1 \leq p \leq \infty$, where $C>0$ is a constant. The above bounds yield
\begin{align}
\| \uk_1\|_{L^p_{x_1,x_2}(\HH)} 
&\leq C \|\UE(t)\|_{L^p_{x_1}} \left( (4 \nu t)^{1/(2p)} + C_{\eta} (4 \nu t)^{1/2} \right)  \leq \CE (\nu t)^{1/(2p)} \label{eq:uk:1:Lp}\\
\|\partial_{1} \uk_1\|_{L^p_{x_1,x_2}(\HH)} 
&\leq C_{\eta}\|\partial_1 \UE(t)\|_{L^p_{x_1}}(\nu t)^{1/(2p)}  \leq \CE (\nu t)^{1/(2p)}  \label{eq:uk:1:d1:Lp} \\
\|\partial_{2} \uk_1\|_{L^p_{x_1,x_2}(\HH)} 
&\leq C_{\eta}\| \UE(t)\|_{L^p_{x_1}}(\nu t)^{1/(2p)-1/2} \leq \CE (\nu t)^{1/(2p)-1/2} \label{eq:uk:1:d2:Lp} \\
\|\partial_{12} \uk_1\|_{L^p_{x_1,x_2}(\HH)} &\leq C_{\eta}\|\partial_1\UE(t)\|_{L^p_{x_1}}(\nu t)^{1/(2p)-1/2}  \leq \CE (\nu t)^{1/(2p)-1/2}\label{eq:uk:1:d12:Lp}
\end{align}
for all $ 1 \leq p \leq \infty$, where $\CE >0$ is a constant that depends on the Euler flow, on $p$, the cutoff function $\eta$, through the constant $C_\eta$, on $\nu_0$ and $T$. We emphasize however only the dependence on the Euler flow.

We moreover have that
\begin{align*}
\partial_t \uk_1 - \nu \Delta \uk_1 
&= - \left( \partial_t \UE(x_1,t) - \nu \partial_{11} \UE(x_1,t) \right) \left( \erfc(z(x_2,t)) -  \sqrt{4 \nu t} \eta(x_2)\right) \notag\\
& \quad + \UE(x_1,t) (\partial_t - \nu \partial_{22}) \left(  \sqrt{4 \nu t}\, \eta(x_2) \right)
\end{align*}
and thus
\begin{align}
\| \partial_t \uk_1 - \nu \Delta \uk_1 \|_{L^2} 
&\leq C_\eta \left( \|\partial_t \UE\|_{L^2} + \nu \|\partial_{11} \UE\|_{L^2} \right) (\nu t)^{1/4} + C_\eta \|\UE\|_{L^2} \nu^{1/2} t^{-1/2}\notag\\
&\leq \CE \left( (\nu t)^{1/4} + \nu^{1/2} t^{-1/2} \right)
 \label{eq:uk:1:heat:L2}
\end{align}
where as before the dependence of all constants on $\nu_0$ and $T$ is ignored.

\subsection{The normal component of the lift $\uk$}
Combining \eqref{eq:uk:incompressible} with \eqref{eq:uk:1:def}, we arrive at 
\begin{align}
\uk_2(x_1,x_2,t) 
&= \partial_1 \UE(x_1,t) \left( \int_0^{x_2}   \erfc(z(y,t)) dy - \sqrt{4 \nu t} \int_0^{x_2} \eta(y) dy \right)\notag\\
&= \sqrt{4 \nu t}\; \partial_1 \UE(x_1,t) \left( \int_0^{z(x_2,t)}  \erfc(z)   dz - \int_0^{x_2} \eta(y) dy  \right) \notag\\
&=: \sqrt{4 \nu t}\; \partial_1 \UE(x_1,t) R(x_2,t).
\label{eq:uk:2:def}
\end{align}
An explicit calculation shows that 
\begin{align*}
R(x_2,t)  = \left(\frac{1}{\sqrt{\pi}} - \int_1^{x_2} \eta(y) dy\right) - \frac{1}{\sqrt{\pi}} \exp\left(-z(x_2,t)^2\right) + z(x_2,t)   \erfc(z(x_2,t)).
\end{align*}
Moreover, note that in view of the choice of $\eta$ in \eqref{eq:eta}, the first term on the right side of the above is identically vanishing for all $x_2 \geq 2$.
It is clear that $R$ obeys
\begin{align*}
R(0,t) = 0 = \lim_{x_2 \to \infty} R(x_2,t),
\end{align*}
and thus we may hope that $R$ is integrable with respect to $x_2$, which is indeed the case.
To see this, first we note that 
\begin{align*}
\|R(t)\|_{L^\infty_{x_2}} \leq \frac{1}{\sqrt{\pi}}.
\end{align*}
Then, we have that 
\begin{align*}
\|R(t)\|_{L^1_{x_2}} 
&\leq \int_0^2 \left|\frac{1}{\sqrt{\pi}} - \int_1^{x_2} \eta(y) dy\right| dx_2 + \frac{1}{\sqrt{\pi}} \int_0^\infty \exp( - z(x_2,t)^2) + z(x_2,t) \erfc(z(x_2,t)) dx_2  \notag\\
&\leq C_\eta + \frac{\sqrt{4\nu t}}{\sqrt{\pi}} \int_0^\infty \exp( - z^2) + z (1 - \erf(z)) dz\notag\\
&\leq C_\eta
\end{align*}
where the dependence of all constants on $\nu_0$ and $T$ is ignored. By interpolation it then follows that 
\begin{align}
\|R(t)\|_{L^p_{x_2}} \leq C_{\eta}
\label{eq:R:Lp}
\end{align}
for all $1 \leq p \leq \infty$. In view of \eqref{eq:R:Lp} and \eqref{eq:uk:1:def}, we have that the bounds
\begin{align}
\|\uk_2\|_{L^p_{x_1,x_2}(\HH)} \leq C_{\eta} \sqrt{4 \nu t} \|\partial_1 \UE\|_{L^p_{x_1}} \leq \CE (\nu t)^{1/2} \label{eq:uk:2:Lp}\\
\|\partial_1 \uk_2\|_{L^p_{x_1,x_2}(\HH)} \leq C_{\eta} \sqrt{4 \nu t} \|\partial_{11} \UE\|_{L^p_{x_1}} \leq \CE (\nu t)^{1/2}\label{eq:uk:2:d1:Lp}
\end{align}
hold for   $1\leq p \leq \infty$, where we have as before suppressed the dependence on $C_\eta$ and $p$ of the constant $C_E$.

Lastly, we obtain from \eqref{eq:uk:incompressible} and \eqref{eq:uk:2:def} that 
\begin{align*}
&(\partial_t - \nu \Delta) \uk_2(x_1,x_2,t) \notag\\
&\quad = \nu \partial_{12} \uk_1(x_1,x_2,t) - \nu \sqrt{4 \nu t}\; \partial_{111} \UE(x_1,t) R(x_2,t) \notag\\
&\qquad + \nu^{1/2} t^{-1/2}  \partial_1 \UE(x_1,t) R(x_2,t)   + \sqrt{4 \nu t}\; \partial_1 \UE(x_1,t) \partial_t R(x_2,t)\notag\\
&\quad = \nu \partial_{12} \uk_1(x_1,x_2,t) - \nu \sqrt{4 \nu t}\; \partial_{111} \UE(x_1,t) R(x_2,t) \notag\\
&\qquad + \nu^{1/2} t^{-1/2}  \partial_1 \UE(x_1,t) R(x_2,t)  - \nu^{1/2} t^{-1/2}\; \partial_1 \UE(x_1,t) z(x_2,t) \erfc(z(x_2,t))
\end{align*}
where we have used that 
\begin{align*}
\partial_t R(x_2,t) = - \frac{1}{2t} z(x_2,t)  \erfc(z(x_2,t)).
\end{align*}
Using \eqref{eq:uk:1:d12:Lp} and \eqref{eq:R:Lp} we conclude that
\begin{align}
\| (\partial_t - \nu \Delta) \uk_2\|_{L^2_{x_1,x_2}(\HH)} 
&\leq C_\eta \nu^{1/2} t^{-1/2} (\nu t)^{1/4} \|\partial_1\UE(t)\|_{L^2_{x_1}}   \notag\\
&\qquad + C_\eta \nu (\nu t)^{1/2} \|\partial_{111} \UE\|_{L^2_{x_1}} + C_\eta \nu^{1/2} t^{-1/2} \|\partial_1 \UE\|_{L^2_{x_1}} \notag\\
&\leq \CE \left( \nu^{1/2} t^{-1/2} + (\nu t)^{1/2} \right)
\label{eq:uk:2:heat:L2}
\end{align}
holds.

\section{Conclusion of the Proof of Theorem~\ref{thm:main}}
\label{sec:proof:2}
Having constructed the corrector function $\uk$, we estimate the terms on the right side of \eqref{eq:v:2}.

\subsection{Bounds for $T_1$, $T_2$, $T_3$,  $T_4$, and $T_5$}
Using \eqref{eq:uk:1:heat:L2} and \eqref{eq:uk:2:heat:L2} we arrive at
\begin{align}
|T_1|
&\leq \|v\|_{L^2} \| (\partial_t - \nu \Delta) \uk\|_{L^2} \notag\\
&\leq \CE \|v\|_{L^2} \left( \nu^{1/2} t^{-1/2} + (\nu t)^{1/4} \right).
\label{eq:T1:bound}
\end{align}
In order to bound $T_2$ we first estimate
\begin{align*}
|T_2|
&\leq  \|\nabla \ue\|_{L^\infty} \| \uk\|_{L^2} \|\uns\|_{L^2}   \notag\\
&\leq  \|\nabla \ue\|_{L^\infty} \| \uk\|_{L^2}  \|\uns_0\|_{L^2} 
\end{align*}
where we have used the $L^2$ energy inequality for the Navier-Stokes solution.  Combining the above with \eqref{eq:uk:1:Lp} and \eqref{eq:uk:2:Lp} we arrive at
\begin{align}
|T_2| \leq  \CE (\nu t)^{1/4} 
\label{eq:T2:bound}
\end{align}
since $\| \uns_0\|_{L^2} \leq C (\|\ue_0\|_{L^2}+1)$, for all $\nu \leq \nu_0$, as we assume $\| \uns_0 - \ue_0\|_{L^2} \to 0$ as $\nu \to 0$.
Similarly to $T_2$, we may estimate
\begin{align}
|T_3|
&\leq  \|\nabla \ue\|_{L^\infty} \| \uk\|_{L^2}   \|v\|_{L^2}  \notag\\
&\leq \CE (\nu t)^{1/2} \|v\|_{L^2}.
\label{eq:T3:bound}
\end{align}
Then, similarly to $T_2$ we estimate $T_4$. We appeal to the energy inequality for the Navier-Stokes solution and estimate \eqref{eq:uk:2:d1:Lp}, which is valid also for $p=\infty$, to conclude that
\begin{align}
|T_4| 
&\leq \| \uns\|_{L^2}^2 \|\partial_1 \uk_2\|_{L^\infty} \notag\\
&\leq \| \uns_0\|_{L^2}^2 \|\partial_1 \uk_2\|_{L^\infty} \notag\\
&\leq \CE (\nu t)^{1/2}.
\label{eq:T4:bound}
\end{align}
\blu{Finally, we estimate $T_5$ by writing $\uns = v + (\ue + \uk)$ so that by the triangle inequality
\begin{align}
|T_5| &\leq C \| v\|_{L^2}^2 \|\partial_1 \uk_1\|_{L^\infty} + C \|v\|_{L^2} \|\ue + \uk\|_{L^\infty} \|\partial_1 \uk_1\|_{L^2} + C \|\ue + \uk\|_{L^\infty}^2 \|\partial_1 \uk_1\|_{L^1}\notag\\
&\leq \CE \|v\|_{L^2}^2 + \CE \|v\|_{L^2} (\nu t)^{1/4} + \CE (\nu t)^{1/2} \notag\\
&\leq \CE \|v\|_{L^2}^2 + \CE (\nu t)^{1/2}.
\label{eq:T5:bound}
\end{align}
In the above estimate we have implicitly used the bounds \eqref{eq:uk:1:Lp}--\eqref{eq:uk:1:d1:Lp} and \eqref{eq:uk:2:Lp}--\eqref{eq:uk:2:d1:Lp}.}

\subsection{Bound for $T_6$}
First we note that by the definition of $\uk_1$ in \eqref{eq:uk:1:def} we have
\begin{align}
|T_6|
&\leq (4 \nu t)^{1/2} \left|\int_{\HH} \uns_1(x_1,x_2,t) \uns_2(x_1,x_2,t) \UE(x_1,t) \eta'(x_2) dx_1 dx_2 \right| \notag\\
& \qquad + \left|\int_{\HH} \uns_1(x_1,x_2,t) \uns_2(x_1,x_2,t) \UE(x_1,t) \partial_{x_2} \erfc(z(x_2,t)) dx_1 dx_2 \right|\notag \\
&\leq \CE (\nu t)^{1/2} \|\uns\|_{L^2}^2 + |T_{6,\nu}| \notag\\
&\leq \CE (\nu t)^{1/2} + |T_{6, \nu}|
\label{eq:T6:1}
\end{align}
where we have used the energy inequality  $\| \uns\|_{L^2} \leq \|\uns_0\|_{L^2} \leq C (1 + \|\ue_0\|_{L^2})$ , and have denoted
\begin{align}
T_{6,\nu} 
&= \int_{\HH} \uns_1(x_1,x_2,t) \uns_2(x_1,x_2,t) \UE(x_1,t) \partial_{x_2} \erfc(z(x_2,t)) dx_1 dx_2 \notag\\
&= - \frac{1}{\sqrt{\pi \nu t}} \int_{\HH} \uns_1(x_1,x_2,t) \uns_2(x_1,x_2,t) \UE(x_1,t) \exp\left(- \frac{x_2^2}{4 \nu t}\right) dx_1 dx_2 \notag\\
&= - \frac{2}{\sqrt{\pi}} \int_{\HH} \uns_1(x_1, \sqrt{4 \nu t} y, t) \uns_2(x_1,\sqrt{4 \nu t} y,t) \UE(x_1,t) \exp(-y^2) dx_1 dy.
\label{eq:T6:2}
\end{align}
The goal is now to show that assumptions \eqref{eq:assume:1}--\eqref{eq:assume:2} imply
\begin{align}
\lim_{\nu \to 0} \int_0^T |T_{6,\nu}(t)| dt = 0
\label{eq:T6:4}
\end{align}
which yields the desired $T_6$ estimate. 

In order to prove \eqref{eq:T6:4}, we fix an $\eps>0$, arbitrary, which in turn fixes a $\rho = \rho(\eps)>0$ such that \eqref{eq:assume:2b} holds. 
\blu{
We then have
\begin{align}
\int_0^T |T_{6,\nu}(t)| dt  
&\leq \CE  \int_0^T\!\!\! \int_{x_2 \geq \rho} \left| \uns_1(x_1,x_2,t) \uns_2(x_1,x_2,t)\right|  \frac{\exp\left( \frac{- x_2^2}{4 \nu t}\right)}{\sqrt{\pi \nu t}} dx_1 dx_2 dt\notag\\
&\quad + \CE \int_0^T\!\!\! \int_{x_2 \leq \rho}  \left| \uns_1(x_1,x_2,t) \uns_2(x_1,x_2,t) \right| \frac{\exp\left(- \frac{x_2^2}{4 \nu t}\right)}{\sqrt{\pi \nu t}} dx_1 dx_2 dt \notag\\
&\leq  \CE \int_0^T \frac{\exp\left(\frac{- \rho^2}{4 \nu t}\right)  }{\sqrt{\pi \nu t}} \|\uns(t)\|_{L^2}^2 dt  + \eps \CE \int_0^T\!\!\! \int_{x_2 \leq \rho} \gamma(x_1,t) \frac{\exp\left(\frac{- x_2^2}{4 \nu t}\right)}{\sqrt{4 \nu t}} dx_1 dx_2 dt \notag\\
&\leq \CE \|\uns_0\|_{L^2}^2  \frac{T \exp\left(\frac{- \rho^2}{4 \nu T}\right)}{\sqrt{\pi \nu T}}  + \eps \CE  \|\gamma\|_{L_1(0,T;L^1_{x_1})} \int_{y \leq \rho/\sqrt{4\nu t}} \exp(-y^2) dy \notag\\
&\leq \CE  \frac{\exp\left(\frac{-\rho^2}{4 \nu T}\right)}{\sqrt{\pi \nu T}}  + \eps \CE  \|\gamma\|_{L^1(0,T;L^1_{x_1})}.
\label{eq:T6:5}
\end{align}
}
By passing $\nu \to 0$ in \eqref{eq:T6:5}, for $\rho$ and $T$ are fixed, we arrive at
\begin{align}
\lim_{\nu \to 0} \int_0^T |T_{6,\nu}(t)| dt \leq \eps \CE  \|\gamma\|_{L^1(0,T;L^1_{x_1})}.
\label{eq:T6:6}
\end{align}
Since $\gamma$ is independent of $\eps$, and $\eps>0$ is arbitrary, \eqref{eq:T6:6} implies \eqref{eq:T6:4} as desired.

\subsection{Proof of Theorem~\ref{thm:main}}
From \eqref{eq:v:2} and \eqref{eq:T1:bound}--\eqref{eq:T6:1} we conclude that 
\begin{align}
\frac{d}{dt} \|v\|_{L^2}^2
&\leq \CE  \|v\|_{L^2}^2 
+ \CE  \nu^{1/2} t^{-1/2}   \|v\|_{L^2} + \CE (\nu t)^{1/4}  + T_{6,\nu}
\label{eq:v:3}
\end{align}
where as usual $\CE$ implicitly depends on $\nu_0$ and $T$.
Upon integrating \eqref{eq:v:3} in time, using \eqref{eq:IC:cond} and \eqref{eq:T6:4} we arrive at 
\begin{align*}
\lim_{\nu \to 0} \|v\|_{L^\infty(0,T;L^2(\HH))} = 0.
\end{align*}
The above yields the proof of \eqref{eq:inviscid:limit}  once we recall that $\uns - \ue = v + \uk$, and that cf.~\eqref{eq:uk:1:Lp} and \eqref{eq:uk:2:Lp} we have $\lim_{\nu \to 0} \|\uk\|_{L^\infty(0,T;L^2(\HH))} = 0$.

\section{Proof of Theorem~\ref{cor:2}}
\label{sec:proof:3}
The proof follows from the proof of Theorem~\ref{thm:main}, as soon as we manage to establish the limit \eqref{eq:T6:4} for the $T_6$ term.
Recall that
\blu{
\begin{align*}
\frac{\sqrt{\pi}}{2} |T_{6,\nu}(t)| =  \left| \int_{\HH} \uns_1(x_1, x_2, t) \uns_2(x_1, x_2,t) \UE(x_1,t) \frac{\exp\left(- \frac{x_2^2}{4 \nu t}\right)}{\sqrt{4 \nu t}} dx_1 dx_2 \right|.
\end{align*}
Fix $\eps>0$ and fix $\rho(\eps)>0$ such that \eqref{eq:assume:G3} holds. Then, as in \eqref{eq:T6:5}, the contribution to $\int_0^T |T_{6,\nu}(t)| dt$ from $\{ x_2 \geq \rho\}$ may be estimated as
\begin{align}
\CE \frac{\exp\left(- \frac{\rho^2}{4 \nu T}\right)}{\sqrt{\pi \nu T}}  \to 0 \quad \mbox{as} \quad \nu \to 0,
\label{eq:VV:1}
\end{align}
for $\eps$, and thus $\rho$, fixed.
}

\blu{
On the other hand, by appealing to both \eqref{eq:assume:1} and \eqref{eq:assume:G3}, the divergence free condition and the fact that $u_2(x_1,0,t)=0$,  the contribution to $\int_0^T |T_{6,\nu}(t)| dt$ from $\{x_2 \leq \rho\}$ may be bounded by
\begin{align}
&\CE \int_0^T \|\uns_1(t) \,{\bf 1}_{0<x_2<\rho}\|_{L^\infty} \int_{\HH} \left| \int_0^{x_2} \partial_1 \uns_1(x_1,y,t) dy \right| dx_1  \frac{\exp\left(- \frac{x_2^2}{4 \nu t}\right)}{\sqrt{4 \nu t}} dx_2 dt \notag\\
&\quad \leq \CE \int_0^T \|\uns_1(t) \, {\bf 1}_{0<x_2<\rho}\|_{L^\infty} \int_{\HH} \int_0^{\rho}  \left|\partial_1 \uns_1(x_1,y,t) \right| dy dx_1  \frac{\exp\left(- \frac{x_2^2}{4 \nu t}\right)}{\sqrt{4 \nu t}} dx_2 dt\notag\\
&\quad \leq \CE \|\uns_1 \, {\bf 1}_{0<x_2<\rho} \|_{L^2(0,T;L^\infty)} \| \partial_1 \uns_1 \, {\bf 1}_{0<x_2<\rho}\|_{L^2(0,T;L^1)}
\notag\\
&\quad \leq \eps \CE \CNS \nu_0.
\label{eq:VV:2}
\end{align}
}

Therefore, by adding the estimates \eqref{eq:VV:1} and \eqref{eq:VV:2}, we obtain
\begin{align*}
\lim_{\nu \to 0} \int_0^T |T_{6,\nu}(t)| dt \leq C \eps
\end{align*}
for any $\eps>0$, as desired.

\blu{
Note that Remark~\ref{rem:L2:x2} follows once we distribute the $x_1$ integrability in \eqref{eq:VV:2}, equally between $\uns_1$ and $\partial_1 \uns_1$ via the Cauchy-Schwartz inequality.
}

\section{Proof of Theorem~\ref{thm:1}}
\label{sec:proof:4}
The proof follows closely that of Theorem~\ref{thm:main}. To avoid redundancy, here we only point out the  main differences. Moreover, for the sake of simplicity, we first consider the case $\delta(\nu t) = \sqrt{4 \nu t}$, which clearly obeys condition \eqref{eq:Wang:assumption} \blu{and the required conditions on $\delta(\nu t)$ stated in the theorem}.

We need to show that \eqref{eq:assume:4} and \eqref{eq:assume:3} imply 
that \eqref{eq:T6:4} holds. Once this is proven, the theorem follows with the same argument as Theorem~\ref{thm:main}. 

\blu{
For this purpose, we may decompose the integral defining $T_{6,\nu}$ into $x_2 > \rho = \sqrt{4 \nu t} (\log(1/\nu))^{1/2}$ and $x_2 < \rho= \sqrt{4 \nu t} (\log(1/\nu))^{1/2}$. The first integral, for $x_2$ away from $0$ is bounded as in \eqref{eq:T6:5} by
\begin{align*}
\CE \| \uns_0\|_{L^2}^2 \int_0^T \frac{ \exp\left( - \frac{\rho^2}{4 \nu t} \right)}{\sqrt{4 \nu t}} dt 
&= \CE \| \uns_0\|_{L^2}^2 \int_0^T \frac{ \exp\left( - \log (1/\nu) \right)}{\sqrt{4 \nu t}} dt  \leq \CE \nu^{1/2} \to 0 
\end{align*}
as $\nu \to 0$. For the contribution to $T_{6,\nu}$ from $x_2 < \rho$, we change variables $y = x_2/\sqrt{4 \nu t}$, so that it remains to show that
\begin{align}
\lim_{\nu \to 0} \int_0^T \!\!\! \int_{y < (\log(1/\nu))^{1/2}} \left| \uns_2(x_1, \sqrt{4 \nu t}\, y ,t)\uns_1(x_1,\sqrt{4 \nu t} \,y ,t)  \UE(x_1,t)  \exp\left(-y^2\right) \right| dx_1 dy dt = 0.
\label{eq:T6:7}
\end{align}
}
Let $M(t) \geq 0$ be defined by
\begin{align*}
\sup_{\nu \in (0,1]}   \| \uns_1 (t)  \uns_2(t) \|_{ L^\infty (\{x_2 \leq \delta(\nu t) (\log(1/\nu))^{1/2}\}) }  = M^2(t).
\end{align*}
By assumption~\eqref{eq:assume:3} we have that $\int_0^T M^2(t) dt < \infty$, and thus the  function 
\begin{align*}
A(x_1,y,t) =  M^2(t) |\UE(x_1,t)| \exp(-y^2) 
\end{align*}
is independent of $\nu$, obeys 
\[
A  \in L^1(dt dy dx_2),
\] 
since the Euler trace $\UE$ is bounded in $L^\infty(0,T;L^1_{x_1}(\RR))$,
and we have that
\begin{align*}
\left|   \uns_2(x_1,  \sqrt{4 \nu t}\, y,t)\uns_1(x_1,  \sqrt{4 \nu t}\, y,t)  \UE(x_1,t)  \exp\left(-y^2\right) \right| \leq A(x_1,y,t)
\end{align*}
for a.e. $(x_1,y,t)$, and all $\nu \in (0,\nu_0]$, in view of assumption \eqref{eq:assume:3}. Thus, in view of \eqref{eq:assume:4}, which guarantees that
\begin{align*}
\lim_{\nu \to 0} \left| \uns_1(x_1,\delta(\nu t)\, y, t) \uns_2(x_1,\delta(\nu t)\, y, t) \UE(x_1,t)  \exp\left(-y^2\right)\right|  = 0
\end{align*} 
we may apply the Dominated Convergence Theorem and conclude that \eqref{eq:T6:7} holds. This concludes the proof of the theorem when $\delta(\nu t) = \sqrt{4 \nu t}$.
 
 To treat the more general case $\delta(\nu t)$ which obeys \eqref{eq:Wang:assumption}, we need to define a corrector \blu{with a different length scale}. For this purpose, \blu{we either recall the construction in~\cite{ConstantinKukavicaVicol15}, or note that we may consider the corrector $\uk$ constructed in \eqref{eq:uk:1:def}--\eqref{eq:uk:2:def},  in which we replace $\sqrt{4 \nu t}$ with $\delta(\nu t)$, and $z(x_2,t)$ with  $x_2/\delta(\nu t)$.
Similar estimates to those in Section~\ref{sec:corrector} show that the} bounds
\begin{align*}
\|\uk\|_{L^p(\HH)} + \blu{t \| \partial_t \uk\|_{L^p(\HH)}} +  \|\partial_1 \uk\|_{L^p(\HH)}  + \|\partial_{11} \uk\|_{L^p(\HH)}  &\leq \CE \delta(\nu t)^{1/p} \\
\| \partial_2 \uk_1\|_{L^p(\HH)} &\leq \CE \delta(\nu t)^{-1 + 1/p} \\
\| \partial_1 \uk_2\|_{L^p(\HH)} &\leq \CE \delta(\nu t)
\end{align*}
hold. \blu{The assumption on $\nu t \delta'(\nu t) / \delta(\nu t) \leq C(T)$ is only needed to ensure that the $\partial_t \uk$ estimate holds.
Although this corrector does not solve the heat equation we may appeal to the usual integration by parts and Poincar\'e inequality trick of~\cite{Kato84b} to bound the $T_1$ term.} It then follows that the terms $T_1,\ldots,T_5$ defined in \eqref{eq:T1:def}--\eqref{eq:T5:def} obey the estimates
\begin{align}
|T_1| &\leq \CE \blu{\frac{\delta(\nu t)^{1/2}}{t}} \|v\|_{L^2} + \CE \nu \delta(\nu t)^{1/2} \|v\|_{L^2} + \frac{\nu}{2} \|\partial_2 v\|_{L^2}^2 + \nu \|v\|_{L^2}^2 + \CE \frac{\nu}{\delta(\nu t)} \label{eq:T1:new} \\
|T_2| &\leq \CE \delta(\nu t)^{1/2} \\
|T_3| &\leq \CE \delta(\nu t)^{1/2} \|v\|_{L^2} \\
|T_4| &\leq \CE \delta(\nu t) \\
|T_5| &\leq \blu{\CE \|v\|_{L^2}^2 + \CE \delta(\nu t)^{1/2}}.
\end{align}
\blu{If there were no $T_6$ term, the proof is completed upon integrating the above bounds on $[0,T]$ and passing $\nu \to 0$, since condition \eqref{eq:Wang:assumption} ensures that the contribution from the first and last terms on the right side of \eqref{eq:T1:new} vanishes in the limit.} For the term $T_6$ we proceed as above, by appealing to \eqref{eq:assume:4}--\eqref{eq:assume:4} and the Dominated Convergence Theorem.  We omit further details.

\section*{Acknowledgements}
\blu{The authors would like to thank Gung-Min Gie and Jim Kelliher for many useful comments about the paper, in particular for pointing out how to remove one unnecessary assumption from an earlier version of Theorem 1.1.}
The work of PC was partially supported by NSF grant DMS-1209394. The work of VV was partially supported by NSF grant DMS-1514771 and by an Alfred P. Sloan Research Fellowship.

\newcommand{\etalchar}[1]{$^{#1}$}

\end{document}